\title{ Ramsey numbers and the size of graphs}
\author{
Benny Sudakov \thanks{ Department of Mathematics, Princeton University, Princeton, NJ 08544, and
Institute for Advanced Study, Princeton. E-mail:
bsudakov@math.princeton.edu.
Research supported in part by NSF CAREER award DMS-0546523, NSF grant
DMS-0355497, USA-Israeli BSF grant, Alfred P. Sloan fellowship, and
the State of New Jersey.
}
}

\date{}
\documentclass [11pt]{article}
\usepackage{graphics}
\usepackage{amsfonts}
\usepackage{amsmath}
\newtheorem{theo}{Theorem}[section]
\newtheorem{prop}[theo]{Proposition}

\begin{document}
\maketitle

\begin{abstract}
For two graph $H$ and $G$, the Ramsey number $r(H, G)$ is the smallest positive integer $n$
such that every red-blue edge coloring of the complete graph $K_n$ on $n$ vertices contains either 
a red copy of $H$ or a blue copy of $G$. Motivated by questions of Erd\H{o}s and Harary, in this note 
we study how the Ramsey number $r(K_s, G)$ depends on the size of the graph $G$. For $s \geq 3$,
we prove that for every $G$ with $m$ edges, $r(K_s,G) \geq c\big(m/\log 
m\big)^{\frac{s+1}{s+3}}$ 
for some positive constant $c$ depending only on $s$. This lower bound improves an earlier result of 
Erd\H{o}s, Faudree, Rousseau, and Schelp, and is tight up to a polylogarithmic factor when $s=3$. 
We also study the maximum value of $r(K_s,G)$ as a function of $m$.
\end{abstract}

\section{Introduction}
For two graphs $H$ and $G$, the {\em Ramsey number} $r(H,G)$ is the smallest positive
integer $n$ such that any red-blue coloring of the edges of the complete
graph $K_n$ on $n$ vertices contains either a red copy of $G$ or a blue copy of
$H$. If $H=G$ we usually denote $r(G, G)$ by $r(G)$. 
The problem of determining or accurately estimating
Ramsey numbers is one of the central problems in modern Combinatorics, and
it has received a considerable attention, see, e.g., \cite{GRS}, \cite{ChangG}. 
In most cases the Ramsey number is estimated in terms of the order (number of vertices) 
of the the graph. However, in the early 80's Erd\H{o}s and Harary  asked about 
the relation between $r(H, G)$ and the sizes (number of edges) 
of the graphs $H$ and $G$. 

The first partial answers to this general problem were obtained by 
Erd\H{o}s, Faudree, Rousseau, and Schelp \cite{EFRS}. They determined up to a constant factor the 
minimum value of $r(G)$ for all graphs of size $m$ and showed that the order of 
magnitude of this minimum is $\Theta(m/\log m)$. They also proved that for fixed $s \geq 3$
there exist constants $c_1, c_2$ such that
$$ c_1 m^{\frac{s}{s+2}} < \min_{e(G)=m}\, r(K_s, G) < c_2\, 
m^{\frac{s-1}{s}}. $$
This estimate is not sharp and there is a large gap between the upper and lower bounds
even when the complete graph is a triangle ($s=3$). In this case Erd\H{o}s \cite{Er} conjectured that 
the upper bound $O(m^{2/3})$ is closer to the truth. Our first result improves the bounds from 
\cite{EFRS} and confirms this conjecture.

\begin{theo}
\label{minimum-ramsey} Let $s\geq 3$ and let $G$ be a graph with
$m$ edges. Then there exists a constant $c$ depending only on $s$
such that
$$r(K_s,G) \geq c\big(m/\log m\big)^{\frac{s+1}{s+3}}.$$ On the other
hand, there exists a graph $G$ of size $m$ such that $r(K_s,G) \leq
O\left( m^{\frac{s-1}{s}}/ \log^{\frac{s-2}{s}}m\right)$.
\end{theo}
In particular, when $s=3$ this determines the minimum value of $r(K_3, G)$ for all graphs $G$ of size
$m$ up to a polylogarithmic factor and shows that
\begin{equation}
\label{s=3}
\Omega\left(\frac{m^{2/3}}{\log^{2/3}m}\right) < \min_{e(G)=m}\, r(K_3, G) < 
O\left(\frac{m^{2/3}}{\log^{1/3}m}\right).
\end{equation}

Another specific question which is part of the general problem mentioned in the first paragraph 
is to bound the maximum value of $r(H, G)$ when the graphs $H$ and $G$ have a given size. 
One of the basic results in Ramsey Theory is the fact
that for the complete graph $G$ with $m$ edges,
$r(G)=2^{\Theta(\sqrt m)}$. A conjecture of Erd\H{o}s \cite{E1} (see also \cite{ChangG}) asserts that
there is an absolute constant $c$ such that for any
graph $G$ with $m$ edges, $r(G) \leq 2^{c \sqrt m}$.
This conjecture is still open. For bipartite graphs it
was proved by Alon, Krivelevich and Sudakov \cite{AKS}. They
also show that for general graphs $G$ with $m$ edges,
$r(G) \leq 2^{c \sqrt m \log m}$ for some absolute positive constant $c$.
For the first off-diagonal case Harary conjectured and Sidorenko \cite{S} proved that 
$r(K_3, G) \leq 2m +1$ for any graph $G$ of size $m$ without isolated vertices.
This inequality is best possible, since $r(K_3, G)= 2m +1$ for any tree with $m$ edges.
Thus for $s=3$ the graphs which maximize $r(K_3, G)$ are very sparse. 
However, for $s>3$ Erd\H{o}s conjectured that exactly the opposite is true and to maximize $r(K_s, G)$ 
over all graphs with $m$ edges one should make $G$ as nearly complete as possible. 
Motivated by this question we obtain the following general upper bound on $r(K_s, G)$ for  graphs $G$ 
of size $m$. 
\begin{theo}
\label{maximum-ramsey} Let $s\geq 3$ and let $G$ be a graph with
$m$ edges and without isolated vertices. Then there exists a constant $c$ depending only on $s$
such that
$$r(K_s,G) \leq c\, m^{\frac{s-1}{2}}/\log^{\frac{s-3}{2}}  m\, .$$
\end{theo}
When $G$ is a clique with $m$ edges it is known by the result of Ajtai,  Koml\'os and 
Szemer\'edi \cite{Aj-KS} that $r(K_s,G) \leq O(m^{\frac{s-1}{2}}/\log^{s-2} m)$. Hence our 
estimate has, up to a polylogarithmic 
factor, similar order of magnitude to the best known upper bound for off-diagonal Ramsey numbers of 
cliques.

The rest of this short note is organized as follows. In the next section we present proofs of our main 
results. The final section contains some concluding remarks and open problems.
Throughout the paper we make no attempts to optimize various
absolute constants. To simplify the presentation, we often omit
floor and ceiling signs whenever these are not crucial.
All logarithms are in the natural base  $e$.

\section{Proofs}

To prove Theorem \ref{minimum-ramsey} we use an approach developed by Krivelevich
\cite{Kr}, which is based on probabilistic arguments together with
large deviation inequalities. The fist inequality we need is a
standard bound of Chernoff (see Appendix A in  \cite{AS}) which
states that if $X$ is a binomially distributed random variable with
parameters $m$ and $p$, then for every $a>0$
$$\mathbb{P}[X-pm<-a] \leq e^{-\frac{a^2}{2pm}}.$$
Another large deviation bound, which we use in the proof, was obtained by Erd\H{o}s
and Tetali \cite{ET} (see also Chapter 8.4 in \cite{AS}).

Let $\Omega$ be a finite set (in our instance it is the set of edges of a complete graph) and 
let $R$ be a random subset of $\Omega$ such that $\mathbb{P}[\omega \in R]=p_\omega$ independently for 
all $\omega \in \Omega$. Let $C_i, i\in I$ be subsets of $\Omega$, where $I$ is some finite index 
set. For every $C_i$ we define $A_i$ to be the event that $C_i \subseteq R$. Let $X_i$ be the 
indicator 
random variable of event $A_i$ and let $X=\sum_{i \in I} X_i$ be the number of 
$C_i \subseteq R$. Finally, let $X_0$ be the maximum number of pairwise disjoint subsets
$C_i$ which belong to $R$. Obviously, $X_0 \leq X$. Let $ \mu$ be the expectation of $X$; then 
Erd\H{o}s and Tetali gave the following bound on the possible size of $X_0$:
$$\mathbb{P}[X_0 \geq k] \leq \frac{\mu^k}{k!} \leq \left(\frac{e\, \mu}{k}\right)^k.$$

\noindent {\bf Proof of Theorem \ref{minimum-ramsey}.}\, Let
$n=\frac{1}{3s^3}\big(m/\log m\big)^{\frac{s+1}{s+3}}$ and consider
coloring the edges of the complete graph $K_n$ such that each edge is 
colored randomly and independently red with probability
$p=\frac{1}{3s}n^{-\frac{2}{s+1}}$ and blue with probability $1-p$.
Let $G_1, \ldots, G_t$ be all subgraph of $K_n$ which are
isomorphic to $G$. The number of such subgraphs $t$ is clearly
bounded by the number of injective functions from $V(G)$ to $K_n$,
which in turn is at most the number of permutations on $n$
elements. Thus $ t \leq n!$. For every subgraph $G_i$, let $X_i$ be
the random variable that counts the number of red edges in $G_i$. By
definition, $X_i$ is binomially distributed with parameters $m$
and $p$. Hence, by Chernoff's inequality
$$\mathbb{P}[X_i<mp/2]=\mathbb{P}[X_i-mp<-mp/2] \leq e^{-mp/8}.$$

Also for every subgraph $G_i$ define $Y_i$ to be the 
number of red cliques of order $s$ which 
share at least one edge with $G_i$. Since $G_i$ has $m$ edges, the number of 
$s$-cliques sharing at least one edge with $G_i$ is bounded by 
$mn^{s-2}$. The probability that an $s$-clique is red is clearly $p^{{s \choose 2}}$.
Therefore 
$$\mathbb{E}[Y_i] \leq mn^{s-2}p^{{s \choose 2}}= mp 
n^{s-2}p^{\frac{(s+1)(s-2)}{2}}=\left(\frac{1}{3s}\right)^{\frac{(s+1)(s-2)}{2}}mp \leq \frac{mp}{9s^2} 
\,.$$ Let $Y'_i$ be the maximum number of edge disjoint red $s$-cliques which share at least one edge 
with $G_i$
Then by the Erd\H{o}s-Tetali inequality we have that 
$$\mathbb{P}\big[Y'_i \geq mp/s^2\big] \leq \left(\frac{e\, \mathbb{E}[Y_i]}{mp/s^2}\right)^{mp/s^2} 
\leq \left(\frac{e}{9}\right)^{mp/s^2} \leq e^{-mp/s^2}.$$
By definition of $n$ and $p$ we have that $mp = (3s^3 n)^{\frac{s+3}{s+1}}p \log m \geq s^2 n \log n 
> 8n \log n$. Therefore the probability that for some index $i$ either $X_i<mp/2$ or $Y'_i \geq mp/s^2$ 
is bounded by $t\big(e^{-mp/8}+e^{-mp/s^2}\big) \leq 2n!\,n^{-n}=o(1)$. In particular there exists a 
red-blue edge coloring of $K_n$ such that for every $1 \leq i \leq t$, subgraph $G_i$ contains at 
least 
$mp/2$ red edges and there are at most $mp/s^2$ edge disjoint red $s$-cliques each sharing at least one 
edge with $G_i$.

Fix such a coloring and let $\Gamma$ be the subgraph of red edges in it. Also let $\cal C$ be the 
maximum (under inclusion) collection of edge disjoint cliques of order $s$ in $\Gamma$. Recolor 
edges in all cliques from $\cal C$ by blue and denote the remaining red graph by $\Gamma'$. Note that by 
recoloring we removed from $\Gamma$ the maximum collection of edge disjoint $s$-cliques. Thus $\Gamma'$ 
contains no clique of order $s$. On the other hand in every subgraph $G_i$ we changed the color of at 
most ${s \choose 2}mp/s^2 <mp/2$ red edges. Since originally $G_i$ had at least $mp/2$ red edges, 
we obtain that every subgraph of $K_n$ isomorphic to $G$ still has at least one red edge. 
This implies that new coloring contains no blue copy of $G$ and no red copy of $K_s$ and completes 
the proof of the first statement of the theorem.

To prove the second part, let $G$ be the union of $\frac{2m}{k(k-1)}$ vertex disjoint cliques of order
$k=m^{1/s}(\log m)^{\frac{s-2}{s}}$. By definition, the number of edges in $G$ is at least $m$. To 
estimate the Ramsey number $r(K_s, G)$ we use the result of 
Ajtai,  Koml\'os and Szemer\'edi \cite{Aj-KS} (see also Theorem 12.17 in \cite{Bol}) which bounds 
off-diagonal Ramsey numbers.
They prove that there exists a constant $c$ such that 
$r(K_s, K_k)\leq c \frac{k^{s-1}}{\log^{s-2} k}$. Let 
$$n=c \frac{k^{s-1}}{\log^{s-2} k}+2m/(k-1)=O\left(m^{\frac{s-1}{s}}/\log^{\frac{s-2}{s}} m\right)$$
and consider any red-blue edge coloring of the complete graph $K_n$. We can assume that there is no red  
$s$-clique, or else we are done. Then, since $n\geq r(K_s, K_k)$, we can find a blue clique of order 
$k$. Delete it from the graph and continue this process. Note that as long as we deleted less than 
$\frac{2m}{k(k-1)}$ cliques of order $k$ the remaining number of vertices is still larger than $r(K_s, 
K_k)$ and we 
can find a new blue $k$-clique. In the end we will find at least $\frac{2m}{k(k-1)}$ blue cliques of 
size $k$, 
i.e., a copy of $G$. This implies that $r(K_s, G) \leq O\left(m^{\frac{s-1}{s}}/\log^{\frac{s-2}{s}} 
m\right)$ and completes the proof.
\hfill $\Box$

\vspace{0.3cm}
\noindent
{\bf Proof of Theorem \ref{maximum-ramsey}.}\, We prove the theorem by induction on $s$. Consider the 
case $s=3$. Although one can use results from \cite{EFRS} and \cite{S} to show that $r(K_3, G) \leq 
O(m)$ we include here the
simple proof that $r(K_3, G) \leq 3m$ for the sake of completeness. Clearly we can assume that $G$ is 
connected, since $r(K_3, G_1 \cup G_2) \leq r(K_3, G_1)+r(K_3, G_2)$. Hence the number of vertices of 
$G$ is at most $m+1$. Let $n=3m$ and suppose that the edges of $K_n$ are red-blue colored with no red
triangle. Pick the vertex with maximum red degree in this coloring and let $X, |X|=t$, be the set of 
its red neighbors. Note that all the edges inside $X$ are blue, since there is no red 
triangle. Partition the vertices of $G$ into two sets $V(G)=V' \cup V''$, where $V'$ consists of the 
$t$ vertices with the highest degree. Since the sum of the degrees in $G$ is $2m$, we have that all 
the 
vertices in $V''$ have degree at most $2m/(t+1)$. Now we will find the blue copy of $G$ as follows.
Embed the vertices of $V'$ into $X$ arbitrarily, and then embed the vertices of $V''$ one by 
one. Given a vertex $v \in V''$, let $Y$ be the set of vertices of $K_n$ where we already embedded 
neighbors of $v$. Since the maximum red degree in the coloring is 
$t$ and $|Y| \leq d(v) \leq 2m/(t+1)$ we have that $K_n$ contains at least $3m -t|Y| \geq m+1-|Y|$ 
vertices which are adjacent to all vertices in $Y$ by blue edges. As the total number of vertices of 
$G$ 
is at most $m+1$, one such vertex is still unoccupied and can be used to embed $v$. Continuing this 
process we find a blue copy of $G$.

Now suppose $s>3$ and by induction we have that $r(K_{s-1}, G) \leq c_1\, 
m^{\frac{s-2}{2}}/\log^{\frac{s-4}{2}}  m$. Let 
$n =c_2 m^{\frac{s-1}{2}}/\log^{\frac{s-3}{2}}  m$, where $c_2$ is a sufficiently large constant 
which depends on $c_1$ and which we fix later. Consider a red-blue coloring of the complete graph 
$K_n$ 
such that there is no red copy of $K_s$. If there is a vertex which has at least
$d=c_1\, m^{\frac{s-2}{2}}/\log^{\frac{s-4}{2}}  m$ red neighbors, then this set cannot  
contain a red copy of $K_{s-1}$. Therefore by the induction hypothesis it will contain a blue copy of 
$G$ and we are done. Thus we can assume that the maximum degree in the red subgraph of $K_n$ is at 
most $d$. Set $k=\sqrt{m\log m}$. It is easy to check that, by definition,
$n=\Omega\left(\frac{k^{s-1}}{\log^{s-2} k}\right)$. Therefore, by choosing $c_2$ large enough and  
using the result of Ajtai,  Koml\'os and Szemer\'edi \cite{Aj-KS} on Ramsey numbers,
we get that $n \geq r(K_s, K_k)$. Hence there exists a set $X$ of $k$ vertices which 
spans only blue edges. Again partition the vertices of $G$ into two sets $V(G)=V' \cup V''$, where 
$V'$ consists of the $k$ vertices with the highest degree. Since the 
sum of the degrees in $G$ is $2m$, we have 
that all the vertices in $V''$ have degree at most $2m/(k+1)$. 
Embed the vertices of $V'$ into $X$ arbitrarily, and then embed the vertices of $V''$ one by
one as follows. Given a vertex $v \in V''$, let 
$Y$ be the set of vertices of $K_n$ where we already embedded
neighbors of $v$. Since the maximum red degree in the coloring is
$d$ and $|Y| \leq d(v) \leq 2m/(k+1)$, by choosing sufficiently large $c_2$, we have that there are at 
least 
$$n-d|Y| \geq n-2md/(k+1) > c_2 m^{\frac{s-1}{2}}/\log^{\frac{s-3}{2}}  m -\frac{2m}{\sqrt{m\log 
m}} \Big(c_1\,m^{\frac{s-2}{2}}/\log^{\frac{s-4}{2}}  m\Big) >2m $$ 
vertices in $K_n$ which are adjacent to all vertices in $Y$ by blue edges. Note that the total number 
of 
vertices of $G$ is at most $2m$, as it has no isolated vertices. Therefore there exists an unoccupied 
vertex of $K_n$ which is connected to all vertices in $Y$ by blue edges. This vertex can be used
to embed $v$. In the end of this procedure we obtain  a blue copy of $G$. This
completes the proof of the theorem. \hfill $\Box$

\section{Concluding remarks}
\begin{itemize}
\item
Let $H$ be a graph with $v_H\geq 3$ vertices and $e_H$ edges. The {\em density} $\rho(H)$ of $H$ is 
defined as $\rho(H)= \frac{e_H-1}{v_H-2}$. Define also 
$$\rho^*(H) =\max_{H' \subseteq H} \rho(H').$$
For example, for the complete graph of order $s$ we have $\rho^*(K_s)=\frac{s+1}{2}$.
The arguments in the proof of Theorem \ref{minimum-ramsey} can be used to obtain the 
following more general result. Since the proof of this statement does not require new ideas
and contains somewhat tedious computations we omit it here.
\begin{theo}
\label{minimum-H-ramsey}
Let $H$ be a fixed graph. Then there exists a constant $c$ depending only on $H$
such that for every graph $G$ with $m$ edges,
$$r(H, G) \geq c\big(m/\log m\big)^{\frac{\rho^*}{1+\rho^*}}.$$
\end{theo}
In addition to the triangle, this result is nearly tight when $H$ is the complete bipartite graph
$K_{p,q}$ with $q \gg p$. Indeed it is easy to check from the definition that if $p$ is fixed and $q 
\rightarrow \infty$ then $\rho^*(K_{p,q}) \rightarrow p$. Therefore for every $p$ and 
$\epsilon>0$ there exists $q$ such that
$\frac{\rho^*(K_{p,q})}{1+\rho^*(K_{p,q})}>\frac{p}{1+p}-\epsilon$. Thus, from Theorem 
\ref{minimum-H-ramsey} we have that 
$r(K_{p,q}, G) \geq \Omega\left(m^{\frac{p}{1+p}-\epsilon}\right)$ for every $G$ with $m$ edges. 
On the other hand, from the result of K\"ovari, S\'os,  Tur\'an
\cite{KST} that $K_{p,q}$-free graphs on $n$ 
vertices can have at most $O(n^{2-1/p})$ edges, it follows that such a graph has an independent set of 
size 
$\Omega\big(n^{1/p}\big)$. This implies that $r(K_{p,q}, K_k) \leq O(k^p)$ (see also \cite{AKS, LZ} 
for 
slightly better estimate). Using this bound together with the argument from the proof of the second 
part of Theorem  \ref{minimum-ramsey}, we can show that if $G$ is the disjoint union of 
$\Theta\big(m^{\frac{p-1}{p+1}}\big)$ cliques of order $m^{1/(p+1)}$ then $r(K_{p,q}, G) \leq 
O\left(m^{\frac{p}{1+p}}\right)$.

\item
For $s=3$ the lower bound in Theorem \ref{minimum-ramsey} is tight up to a multiplicative factor of 
$\log^{1/3} m$. It would be very interesting to close this gap. We think that our upper bound 
in (\ref{s=3}) is closer to the truth and there exists an absolute constant $c$ 
such that $r(K_3, G) \geq c m^{2/3}/\log^{1/3} m$ for every graph $G$ of size $m$.
To prove this one might try to use an approach based on the semi-random method which was developed
by Kim \cite{Kim} to determine the asymptotic behavior of Ramsey numbers $r(K_3, K_k)$.
\item
It would be interesting to extend an upper bound in Theorem \ref{maximum-ramsey} to the general case
when $H$ and $G$ are arbitrary graphs with sizes $t$ and $m$ and with no isolated 
vertices. We conjecture 
that if $t$ is fixed and $m$ is sufficiently large then 
$$r(H, G)\leq m^{O(\sqrt{t})}.$$ 
This estimate if true is tight up to a constant in the exponent, since 
the known lower bounds (see \cite{Sp, Kr}) on 
off-diagonal Ramsey numbers imply that $r(H, G) \geq m^{\Omega(\sqrt{t})}$ when $H$ and $G$ are 
complete graphs with $t$ and $m$ edges respectively. 
In  \cite{AKS} it was proved that if $H$ is a graph with chromatic number $\ell$ and maximum degree $d\geq \ell$ 
then for all sufficiently large $m$, $r(H, K_{2m}) \leq m^{\ell d}$.
Using this estimate it is easy to obtain the following partial result, which shows that our conjecture holds if the
chromatic number of $H$ is a fixed constant.

\begin{prop}
\label{maximum-H-ramsey}
Let $H$ and $G$ be two graphs with no isolated vertices such that the size of $G$ is $m$, the size of 
$H$ is $t$ and $H$ has chromatic number $\ell \geq 2$.
Then there exists a constant $c$ depending only on $\ell$ such that for sufficiently large $m$,
$r(H, G)\leq m^{c\sqrt{t}}$.
\end{prop}

\noindent
{\bf Sketch of proof.}\, We use induction on $t$. Let $v$ be the vertex of maximum degree in $H$.
Since $G$ has $m$ edges, it has at most $2m$ vertices. Therefore, if
the maximum degree of $H$ is at most $2\sqrt{t}$, it follows from the above cited estimate in \cite{AKS}
that $r(H, G)\leq m^{2\sqrt{t}\ell}$.
Otherwise, the degree of $v$ in $H$ is larger than $2\sqrt{t}$. Delete it and denote $H_1=H\setminus\{v\}$. This 
graph has $t_1 \leq t-2\sqrt{t}$ edges  and $\sqrt{t_1} \leq \sqrt{t}-1$. 

Let $n=m^{2\sqrt{t}\ell}$ and consider red-blue
coloring of the edges of the complete graph $K_n$. Since $G$ has at most $2m$ vertices,
we can assume that there is no red $K_{2m}$ in this coloring.
Therefore, by Tur\'an's theorem, there is a vertex $x$ in $K_n$, whose blue degree is at least 
$n/{2m}\gg m^{(2\sqrt{t}-2)\ell} \geq m^{2\sqrt{t_1} \ell}$. Let $U$ be the set of blue neighbors of $x$. 
Clearly, this set contains no blue copy of $H_1$.
Now we can use induction to conclude that it contains a red copy of $G$.
\hfill $\Box$
\end{itemize}


\begin{thebibliography}{99}

\bibitem{Aj-KS}
M. Ajtai,  J. Koml\'os and E. Szemer\'edi,
A note on Ramsey numbers,
{\em J. Combinatorial Theory Ser. A} 29 (1980), 354--360.

\bibitem{AS}
N. Alon and J. Spencer, {\bf The Probabilistic Method}, 2nd ed.,
Wiley, New York, 2000.

\bibitem{AKS}
N. Alon, M. Krivelevich and B. Sudakov, Tur\'an numbers of bipartite
graphs and related Ramsey-type questions , {\em Combinatorics,
Probability and Computing} 12 (2003), 477-494.

\bibitem{Bol}
B. Bollob\'as, {\bf Random Graphs}, 2nd ed., Cambridge Studies in
Advanced Mathematics, 73, Cambridge University Press, Cambridge,
2001.

\bibitem{ChangG}
F. Chung and R. Graham,
{\bf Erd\H{o}s on graphs. His legacy of unsolved
problems}, A K Peters Ltd., Wellesley, MA, 1998.


\bibitem{Er}
P. Erd\H{o}s, Solved and unsolved problems in combinatorics and
combinatorial number theory, Proc. 12th Southeastern Conference on
Combinatorics, Graph Theory and Computing, Vol. I, {\em Congr.
Numer.} 32 (1981), 49--62.

\bibitem{E1} P. Erd\H{o}s,
On some problems in graph theory, combinatorial analysis and
combinatorial number theory, in: {\em Graph theory and combinatorics
(Cambridge, 1983)}, Academic Press, London, 1984, 1--17.


\bibitem{EFRS}
P. Erd\H{o}s, R. Faudree, C. Rousseau and R. Schelp, A Ramsey
problem of Harary on graphs with prescribed size, {\em  Discrete
Math.} 67 (1987), 227--233.

\bibitem{ET}
P. Erd\H{o}s and P. Tetali, Representation of integers as the sum
of $k$ terms, {\em Random Structures and Algorithms} 1 (1990),
245--261.

\bibitem{GRS} R. Graham, B. Rothschild and J. Spencer,
{\bf Ramsey theory}, $2^{nd}$ ed., Wiley, New York, 1990.

\bibitem{Kim}
J.H. Kim, 
The Ramsey number $R(3,t)$ has order of magnitude $t^2/\log t$,
{\em Random Structures Algorithms}  7  (1995), 173--207. 

\bibitem{KST} T. K\"ovari, V.T. S\'os and P. Tur\'an, On a problem
of K. Zarankiewicz,{\em Colloquium Math.} 3 (1954), 50-57.

\bibitem{Kr}
M. Krivelevich, Bounding Ramsey numbers through large deviation
inequalities, {\em Random Structures and Algorithms} 7 (1995),
145--155.

\bibitem{LZ}
Y. Li and W. Zang, 
Ramsey numbers involving large dense graphs and bipartite Tur\'a¡n 
numbers, {\em  J. Combinatorial Theory Ser. B}  87  (2003), 280--288.

\bibitem{S}
A. Sidorenko,
The Ramsey number of an $n$-edge graph versus triangle is at most $2n+1$,
{\em J. Combinatorial Theory Ser. B} 58 (1993), 185--196.

\bibitem{Sp}
J. Spencer, Asymptotic lower bounds for Ramsey numbers,
{\em Discrete Mathematics} 20 (1977), 69--76.

\end{thebibliography}
\end{document}